\documentclass[11pt]{amsart}

\usepackage{amsmath,amssymb,amsthm,mathtools}
\usepackage{graphicx}
\usepackage{float}
\usepackage{tikz}
\usetikzlibrary{shapes.geometric}
\usepackage[colorlinks=true,linkcolor=blue,citecolor=blue,urlcolor=blue]{hyperref}

\newtheorem{proposition}{Proposition}[section]
\newtheorem{theorem}[proposition]{Theorem}
\newtheorem{lemma}[proposition]{Lemma}

\theoremstyle{definition}
\newtheorem{definition}[proposition]{Definition}
\newtheorem{remark}[proposition]{Remark}

\newcommand{\Iperp}{\mathcal{I}_{\perp}}
\newcommand{\Ical}{\mathcal{I}}
\newcommand{\Pcal}{\mathcal{P}}
\newcommand{\RP}{\mathbb{RP}}
\newcommand{\R}{\mathbb{R}}
\newcommand{\Sphere}{\mathbb{S}}
\newcommand{\Torus}{\mathbb{T}}
\newcommand{\diam}{\operatorname{diam}}

\title[New bounds on simplex-avoiding sets in two dimensions]{New bounds on simplex-avoiding sets\\ in two dimensions}

\author{Domonkos Czifra, \'Akos D\'ucz, P\'eter Madarasi,\\ D\'aniel Varga, P\'al Zs\'amboki}

\address{Domonkos Czifra, HUN-REN Alfr\'ed R\'enyi Institute of Mathematics, Re\'altanoda utca 13-15, 1053 Budapest, Hungary}
\email{czifra.domonkos@renyi.hu}

\address{\'Akos D\'ucz, HUN-REN Alfr\'ed R\'enyi Institute of Mathematics, Re\'altanoda utca 13-15, 1053 Budapest, Hungary}
\email{akos@renyi.hu}

\address{P\'eter Madarasi, HUN-REN Alfr\'ed R\'enyi Institute of Mathematics, Re\'altanoda utca 13-15, 1053 Budapest, Hungary; and Department of Operations Research, ELTE E\"otv\"os Lor\'and University, P\'azm\'any P.\ s.\ 1/c, 1117 Budapest, Hungary}
\email{madarasi@renyi.hu}

\address{D\'aniel Varga, HUN-REN Alfr\'ed R\'enyi Institute of Mathematics, Re\'altanoda utca 13-15, 1053 Budapest, Hungary}
\email{daniel@renyi.hu}

\address{P\'al Zs\'amboki, HUN-REN Alfr\'ed R\'enyi Institute of Mathematics, Re\'altanoda utca 13-15, 1053 Budapest, Hungary}
\email{zsamboki.pal@renyi.hu}

\date{\today}

\begin{document}

\begin{abstract}
We study extremal densities of measurable subsets of metric spaces avoiding prescribed finite geometric configurations. We prove two new upper bounds for three-point avoidance problems: a measurable subset of $\Sphere^2$ containing no orthonormal basis has normalized spherical measure at most $7/11$, and a measurable subset of $\R^2$ containing no unit equilateral triangle has upper Banach density at most $9/16$. The first bound appears to be the first improvement over the trivial averaging bound for this problem, while the second improves the previously best known upper bound.

The proofs are applications of the geometric fractional chromatic number method first introduced in~\cite{AmbrusCsiszarikMatolcsiVargaZsamboki2024}. For this purpose, we formulate the method in a framework that applies both to compact homogeneous spaces and to Euclidean spaces with upper Banach density.

\end{abstract}

\maketitle

\section{Introduction}

We study the maximal density of measurable subsets of
geometric spaces avoiding prescribed finite configurations.
Following~\cite{CastroSilva2021}, we refer to such
extremal quantities as independence densities, generalizing
the independence ratio of a graph or hypergraph.

Our first problem concerns subsets of the unit sphere
$\Sphere^2$ containing no orthonormal basis. Let
$m_\perp(\Sphere^2)$ denote the supremum of the normalized
spherical measures of such sets. The elementary averaging
argument gives $m_\perp(\Sphere^2)\le 2/3$.
We prove the following improvement.

\begin{theorem}\label{thm_sphere}
We have
\[
m_\perp(\Sphere^2)\le \frac{7}{11} \approx 0.63636.
\]
\end{theorem}

The problem has previously been considered
in~\cite{Perez2019,
CastroSilvaOliveiraSlotVallentin2022,
Zakharov2023,WilliamsConstantin2025},
mostly in connection with lower bounds or
higher-dimensional variants.
To the best of our knowledge, our bound is the first
explicitly established improvement over $2/3$ for this
independence-density problem. We could deduce a
weaker improvement from the work of Conway and
Kochen~\cite{ConwayKochen2006}: their configuration of
40 orthonormal triples implies, by averaging,
$m_\perp(\Sphere^2)\le 79/120 \approx 0.65833$.

Our second problem concerns measurable subsets of
$\R^2$ containing no unit equilateral triangle.
Let $m_\Delta(\R^2)$ denote the supremum of their
upper Banach densities. A construction of density $1/2$
is conjectured to be optimal. The previously best known
upper bound was $0.64355$, obtained
in~\cite{CastroSilvaOliveiraSlotVallentin2022}.
We improve it as follows.

\begin{theorem}\label{thm_plane}
We have
\[
m_\Delta(\R^2)\le \frac{9}{16} = 0.5625.
\]
\end{theorem}

Both of our results are obtained using the geometric fractional
chromatic number (GFCN) method, introduced
in~\cite{AmbrusCsiszarikMatolcsiVargaZsamboki2024}
and subsequently developed
in~\cite{MatolcsiRuzsaVargaZsamboki2025,
LevMatolcsiPachVarga2025,CzifraEtAl2026}.
The method considers random isometric placements of finite
point configurations, called probes. The resulting
distributions of intersection patterns satisfy linear
constraints arising from the forbidden configurations
and geometric congruences. Solving the corresponding
linear programs yields upper bounds on independence
densities.

Our bounds use a 15-point probe in the projective
plane and a 14-point probe in the Euclidean plane.
The probes were found by computer search, and the bounds
are certified by exact rational dual solutions.
The supplementary material~\cite{SupplementaryMaterial}
contains the certificates and an exact-arithmetic verifier.

Section~\ref{sec:probes} gives an informal overview of the GFCN method. In Sections~\ref{sec:compact} and~\ref{sec:euclidean},
we formulate the method for compact homogeneous
spaces and Euclidean spaces, respectively.
Sections~\ref{sec:sphere} and~\ref{sec:plane}
establish the two bounds above.

\section{The probe method: an informal overview}\label{sec:probes}

The geometric fractional chromatic number (GFCN) method was introduced
in~\cite{AmbrusCsiszarikMatolcsiVargaZsamboki2024}.
Here we describe its
probabilistic, probe-based interpretation, introduced
in~\cite{CzifraEtAl2026}. The basic idea is to study a hypothetical
configuration-avoiding set by testing it with randomly placed finite point
configurations, which we call \emph{probes}. Let $A$ be a Borel set avoiding
the forbidden configurations, and let $X$ be a nonempty finite probe. We
place a congruent copy of $X$ at random and record which points of $X$ land
in $A$. This produces a random subset $Y\subseteq X$, or equivalently a
probability distribution $\gamma_{A,X}$ on the subsets of $X$. The placement
rule is chosen to respect the symmetries of the ambient space. In Euclidean
spaces, we initially work with periodic sets; the precise constructions are
given in Sections~\ref{sec:compact} and~\ref{sec:euclidean}.

There are two sources of constraints on the observation distribution.
The first is combinatorial: an observed subset $Y$ cannot contain a
forbidden configuration, so $\gamma_{A,X}(Y)=0$ for every such $Y$.
Of course, all observation probabilities are nonnegative and sum to one.
These restrictions depend only on which subsets of $X$ contain forbidden
configurations.

The second source is geometric. For any $U\subseteq X$, the probability
that \emph{all} points of $U$ land in $A$, regardless of what happens to
the remaining points of $X$, is
\[
\overline\gamma_{A,X}(U)
  =\sum_{Y:\,U\subseteq Y\subseteq X}\gamma_{A,X}(Y).
\]
If two subsets $U,V\subseteq X$ are congruent under the relevant group
of isometries, their random placements have the same distribution, and
therefore
\[
\overline\gamma_{A,X}(U)=\overline\gamma_{A,X}(V).
\]
These equalities apply to subconfigurations of any size, not just
individual points or forbidden configurations. They use geometric
information about $X$ that is not captured by its forbidden-configuration
hypergraph. The strength of the method comes from combining these
geometric consistency relations with the combinatorial restrictions
on the possible observations.

Every point of the probe has the same probability of landing in $A$,
which equals the density of $A$ (or its periodic density in the Euclidean
setting). We can therefore bound the density by maximizing the common
singleton probability $\overline\gamma(\{x_0\})$ over all probability
distributions on admissible subsets of $X$ satisfying the congruence
equalities. This is a finite linear program. Every avoiding set under
consideration yields a feasible observation distribution, so the LP
optimum bounds its density.
% In Euclidean spaces, a periodic approximation argument extends this conclusion to arbitrary measurable avoiding sets when the forbidden configurations have uniformly bounded diameter.

This density formulation is equivalent to the usual GFCN formulation by
a simple change of normalization.
In the latter, one assigns nonnegative
weights $w(Y)$ to admissible subsets, requires
$\sum_{Y:\,x\in Y}w(Y)=1$ for every $x\in X$, imposes the same congruence
relations, and minimizes the total weight $\sum_Y w(Y)$. If a probability
distribution $\gamma$ has common singleton marginal $p>0$, then
$w(Y)=\gamma(Y)/p$ has total weight $1/p$. Conversely, normalizing such
weights to have total mass one gives a probability distribution with
singleton marginal $1/\sum_Y w(Y)$. Thus the optimum of our density LP
is the reciprocal of the corresponding GFCN.
If we drop the congruence constraints for subsets of size at least
two from the coloring formulation, retaining only the equal singleton marginals, then we recover the
ordinary fractional chromatic number of the hypergraph on $X$
whose edges are the forbidden configurations
(see~\cite{ScheinermanUllman1997}).
Thus the additional geometric consistency constraints are precisely
what distinguishes GFCN from standard fractional coloring. In this paper, we use the probability
normalization because it directly expresses the desired density bound.

The practical challenge is to find probes for which this linear program
gives strong bounds. In our applications, we search for such probes
computationally and verify the resulting bounds using exact rational %TODO2: Ez is volt már az intróban.
certificates. The search procedure is described in
Appendix~\ref{app:search}.

A related use of probabilistic symmetrization appears in the recent proof by OpenAI~\cite{OpenAI2026} that the Euclidean plane is not five-colorable. Their proof averages over translations and rotations to construct invariant probability distributions on colorings, followed by a transfer to weak measurable colorings. Their arguments differ substantially from our finite-probe approach, but both methods exploit geometric symmetries through probabilistic averaging.

\section{Probes in compact spaces}\label{sec:compact}

Let \(S\) be a compact metric space equipped with a prescribed Borel
probability measure \(\delta\).  Let \(G\) be a compact group acting
continuously and transitively on \(S\) by isometries, and assume that \(\delta\)
is \(G\)-invariant.  We write \(\mu\) for normalized Haar measure on \(G\).  The
random placement procedure is to choose \(g\in G\)
according to \(\mu\), and to send each point \(x\in S\) to \(gx\).

\begin{lemma}\label{lem:ambient-measure}
For \(x\in S\), let \(e_x:G\to S\) be the map \(e_x(g)=gx\).  The pushforward
measure \((e_x)_*\mu\) is equal to \(\delta\).
\end{lemma}

\begin{proof}
The measure \((e_x)_*\mu\) is \(G\)-invariant, since left translation preserves
Haar measure.  A compact transitive \(G\)-space has a unique \(G\)-invariant
Borel probability measure: identify \(S\) with \(G/G_x\) and push Haar measure
forward through the quotient map.  Since \(\delta\) is \(G\)-invariant, the two
measures are equal.
\end{proof}

Thus the Haar-random placement rule is compatible with the density measure:
for every Borel set \(A\subseteq S\), \(\delta(A)\) is the probability that a
randomly placed point lands in \(A\).

Let \(\Pcal\) be a \(G\)-invariant family of nonempty finite subsets of \(S\).  Thus, if
\(P\in\Pcal\), then \(gP\in\Pcal\) for every \(g\in G\).  Examples include the
pairs of points at distance \(1\) in a metric space, or the triples forming
orthonormal bases on \(\Sphere^2\).  A set \(A\subseteq S\) is
\(\Pcal\)-independent if no member of \(\Pcal\) is contained in \(A\).  We write
\[
m_{\Pcal}(S)
\coloneq
\sup\{\delta(A) : A\subseteq S \text{ is Borel and }\Pcal\text{-independent}\}.
\]

For a probe \(X\subseteq S\), and for \(Y\subseteq X\), define
\[
E_{A,X}(Y)
\coloneq
\{g\in G : gx\in A \iff x\in Y \text{ for every }x\in X\}
\]
and
\[
\gamma_{A,X}(Y)\coloneq\mu(E_{A,X}(Y)).
\]
This is the probability that a random placement of \(X\) produces exactly the
observed pattern \(Y\).  For a finite set \(Z\subseteq S\), define
\[
\pi_A(Z)\coloneq\mu(\{g\in G : gZ\subseteq A\}).
\]
Since \(g\mapsto gx\) is continuous, these are probabilities of Borel events
whenever \(A\) is Borel.

For \(Z\subseteq X\), define
\[
\overline{\gamma}_{A,X}(Z)
\coloneq
\sum_{Y : Z\subseteq Y\subseteq X}\gamma_{A,X}(Y).
\]
Then
\begin{equation}\label{eq:aggregate-compact}
\pi_A(Z)=\overline{\gamma}_{A,X}(Z).
\end{equation}
Indeed, the event \(gZ\subseteq A\) is the disjoint union of the events
\(E_{A,X}(Y)\) over all \(Y\supseteq Z\).

Let
\[
\Ical_{\Pcal}(X)\coloneq\{Y\subseteq X : Y\text{ is }\Pcal\text{-independent}\}.
\]
If \(A\) is \(\Pcal\)-independent, then \(\gamma_{A,X}(Y)=0\) for every
\(Y\notin \Ical_{\Pcal}(X)\), since otherwise a forbidden configuration inside
\(Y\) would be placed inside \(A\).  Also, if \(U,V\subseteq X\) are congruent,
meaning \(V=hU\) for some \(h\in G\), then
\[
\overline{\gamma}_{A,X}(U)=\overline{\gamma}_{A,X}(V),
\]
because \(\pi_A(V)=\pi_A(U)\) by right invariance of Haar measure.  Finally, for
every \(x\in X\),
\[
\overline{\gamma}_{A,X}(\{x\})=\pi_A(\{x\})=\delta(A).
\]

\begin{definition}
Fix a probe \(X\subseteq S\) and a point \(x_0\in X\).  Let
\(\Gamma_{\Pcal}(X)\) be the optimum of the following linear program.  The
variables are \(\gamma(Y)\), \(Y\in \Ical_{\Pcal}(X)\), and
\[
\overline{\gamma}(Z)\coloneq
\sum_{\substack{Y\in \Ical_{\Pcal}(X)\\ Z\subseteq Y}}\gamma(Y).
\]
Maximize \(\overline{\gamma}(\{x_0\})\), subject to
\[
\gamma(Y)\ge 0,\qquad
\sum_{Y\in \Ical_{\Pcal}(X)}\gamma(Y)=1,
\]
and
\[
\overline{\gamma}(U)=\overline{\gamma}(V)
\qquad
\text{whenever } U,V\subseteq X \text{ are congruent}.
\]
\end{definition}

Since all singleton subsets of $X$ are congruent, the optimum does not depend on the choice of $x_0$.

\begin{theorem}\label{thm:compact-probe-bound}
For every probe \(X\subseteq S\),
\[
m_{\Pcal}(S)\le \Gamma_{\Pcal}(X).
\]
\end{theorem}

\begin{proof}
Every Borel \(\Pcal\)-independent set \(A\) gives a feasible point
\(\gamma_{A,X}\) of the linear program above.  Its objective value is
\[
\overline{\gamma}_{A,X}(\{x_0\})=\delta(A).
\]
Taking the supremum over \(A\) proves the theorem.
\end{proof}

\begin{remark}
If one wants to take the supremum over all \(\delta\)-measurable independent
sets instead of Borel ones, the value is unchanged.  Indeed, \(\delta\) is a
regular Borel probability measure on the compact metric space \(S\), and every
measurable independent set contains compact independent subsets of arbitrarily
close measure.
\end{remark}

\section{Probes in Euclidean spaces}\label{sec:euclidean}

We next record the corresponding form of the method for Euclidean spaces.  Let
\(\operatorname{E}(n)=O(n)\ltimes\R^n\) be the Euclidean isometry group, and
let \(\Pcal\) be an \(\operatorname{E}(n)\)-invariant family of nonempty finite subsets
of \(\R^n\).  Thus, if \(P\in\Pcal\), then \(RP+v\in\Pcal\) for every
\(R\in O(n)\) and \(v\in\R^n\).  A measurable set \(A\subseteq\R^n\) is
\(\Pcal\)-independent if no member of \(\Pcal\) is contained in \(A\).

For a measurable \(A\subseteq\R^n\), define its upper Banach density by
\[
\overline d(A)
\coloneq
\limsup_{T\to\infty}\sup_{u\in\R^n}
\frac{\lambda_n(A\cap(u+[0,T]^n))}{T^n},
\]
where \(\lambda_n\) denotes \(n\)-dimensional Lebesgue measure.  We write
\[
\begin{aligned}
m_{\Pcal}(\R^n)
&\coloneq
\sup\{\,\overline d(A) : A\subseteq\R^n
\text{ is Lebesgue measurable} \\
&\hspace{4.5cm}
\text{and }A\text{ is }\Pcal\text{-independent}\,\}.
\end{aligned}
\]

First consider periodic independent sets.  Let \(L\subseteq\R^n\) be a
full-rank lattice, let
\[
\Torus_L\coloneq\R^n/L,
\]
and let \(q_L:\R^n\to\Torus_L\) be the quotient map.  If
\(A=q_L^{-1}(A_L)\), where \(A_L\subseteq\Torus_L\) is Borel, its periodic
density is
\[
\delta_L(A)\coloneq m_L(A_L),
\]
where \(m_L\) is normalized Haar measure on \(\Torus_L\).

The random placement procedure for an \(L\)-periodic set is to choose
\[
(R,t)\in O(n)\times\Torus_L
\]
according to product Haar measure and to place \(x\in\R^n\) at
\[
\phi_{R,t}(x)\coloneq q_L(Rx)+t\in\Torus_L.
\]
For a probe \(X\subseteq\R^n\), \(Y\subseteq X\), and a periodic Borel
set \(A=q_L^{-1}(A_L)\), define
\[
\gamma_{A,X}(Y)
\coloneq
\Pr\bigl(\phi_{R,t}(x)\in A_L \iff x\in Y
\text{ for every }x\in X\bigr).
\]
For a finite set \(Z\subseteq\R^n\), define
\[
\pi_A(Z)
\coloneq
\Pr\bigl(\phi_{R,t}(z)\in A_L \text{ for every }z\in Z\bigr).
\]
These are probabilities of Borel events, since
\((R,t)\mapsto q_L(Rx)+t\) is continuous.

For \(Z\subseteq X\), define \(\overline{\gamma}_{A,X}(Z)\) by the same sum as in Section~\ref{sec:compact}. As before, \(\pi_A(Z)=\overline{\gamma}_{A,X}(Z)\), since the event that all points of \(Z\) land in \(A_L\) is the disjoint union of the exact observations containing \(Z\).

If \(A\) is \(\Pcal\)-independent, then \(\gamma_{A,X}(Y)=0\) whenever \(Y\notin\Ical_{\Pcal}(X)\), because any placement observing \(Y\) would lift a forbidden configuration \(P\subseteq Y\) to \(RP+\tilde t\subseteq A\) for a lift \(\tilde t\) of \(t\).

If \(U,V\subseteq X\) are Euclidean-congruent, then
\[
\pi_A(U)=\pi_A(V).
\]
Indeed, if \(V=R_0U+v\), then the change of variables
\[
(R,t)\mapsto (RR_0,\ t+q_L(Rv))
\]
preserves product Haar measure on \(O(n)\times\Torus_L\).  Finally,
\[
\pi_A(\{x\})=\delta_L(A)
\]
for every \(x\in\R^n\), because \(q_L(Rx)+t\) is uniformly distributed on
\(\Torus_L\) for every fixed \(R\).

For a Euclidean probe \(X\subseteq\R^n\) and \(x_0\in X\), define \(\Gamma_{\Pcal}(X)\) by the linear program in Section~\ref{sec:compact}, with congruence interpreted as Euclidean congruence. The preceding arguments show that \(\gamma_{A,X}\) satisfies its constraints whenever \(A\) is an \(L\)-periodic \(\Pcal\)-independent set.

\begin{proposition}\label{prop:periodic-euclidean-bound}
If \(A=q_L^{-1}(A_L)\subseteq\R^n\), with \(A_L\subseteq\Torus_L\) Borel, is
\(L\)-periodic and \(\Pcal\)-independent, then
\[
\delta_L(A)\le \Gamma_{\Pcal}(X)
\]
for every probe \(X\subseteq\R^n\).
\end{proposition}

\begin{proof}
The actual observation probabilities \(\gamma_{A,X}\) form a feasible point of
the linear program defining \(\Gamma_{\Pcal}(X)\), and the objective value is
\[
\overline{\gamma}_{A,X}(\{x_0\})=\pi_A(\{x_0\})=\delta_L(A).
\]
\end{proof}

The following lemma gives a periodic approximation when the forbidden configurations have uniformly bounded diameter.

Write $D_{\Pcal}\coloneq\sup_{P\in\Pcal}\diam(P)$, with $D_{\Pcal}=0$ when $\Pcal=\emptyset$.

\begin{lemma}\label{lem:periodic-reduction-euclidean}
Assume \(D_{\Pcal}<\infty\).  Let \(A\subseteq\R^n\) be Lebesgue measurable and
\(\Pcal\)-independent.  For every \(\varepsilon>0\), there are a full-rank
lattice \(L\) and an \(L\)-periodic Borel \(\Pcal\)-independent set
\(B\subseteq\R^n\) such that
\[
\delta_L(B)\ge \overline d(A)-\varepsilon.
\]
\end{lemma}

\begin{proof}
Choose \(\rho>D_{\Pcal}/2\).  Since the boundary-strip proportion tends to
zero as \(T\to\infty\), the definition of \(\overline d(A)\) lets us choose
\(T>2\rho\) and \(u\in\R^n\) such that
\[
\frac{\lambda_n(A\cap(u+[0,T]^n))}{T^n}
\ge \overline d(A)-\frac{\varepsilon}{3}
\]
and such that the boundary strip
\[
(u+[0,T]^n)\setminus(u+[\rho,T-\rho]^n)
\]
has measure at most \(\varepsilon T^n/3\).  By inner regularity, choose a
compact set
\[
K\subseteq A\cap(u+[\rho,T-\rho]^n)
\]
with
\[
\lambda_n(K)\ge
\lambda_n(A\cap(u+[\rho,T-\rho]^n))-\frac{\varepsilon T^n}{3}.
\]
Let \(L=T\mathbb{Z}^n\), and set
\[
B\coloneq K+L.
\]
Then \(B\) is \(L\)-periodic and Borel.  Distinct translates \(K+\ell\),
\(\ell\in L\), are separated by distance at least \(2\rho>D_{\Pcal}\).  Hence
any member of \(\Pcal\) contained in \(B\) would have to lie in a single
translate of \(K\).  Translating it back would give a member of \(\Pcal\)
contained in \(K\subseteq A\), contradicting the \(\Pcal\)-independence of
\(A\).

The quotient map \(q_L\) is injective on \(K\), and \(q_L(K)\) is compact in
\(\Torus_L\).  Therefore \(B=q_L^{-1}(q_L(K))\) and
\[
\delta_L(B)=m_L(q_L(K))=\frac{\lambda_n(K)}{T^n}
\ge \overline d(A)-\varepsilon.
\]
\end{proof}

\begin{theorem}\label{thm:euclidean-probe-bound}
Assume \(D_{\Pcal}<\infty\).  Then, for every probe
\(X\subseteq\R^n\),
\[
m_{\Pcal}(\R^n)\le \Gamma_{\Pcal}(X).
\]
\end{theorem}

\begin{proof}
Apply Lemma~\ref{lem:periodic-reduction-euclidean} to a measurable
\(\Pcal\)-independent set \(A\), then apply
Proposition~\ref{prop:periodic-euclidean-bound} to the resulting periodic set,
and let \(\varepsilon\to0\).
\end{proof}

\section{Orthonormal basis-avoiding sets on the sphere}\label{sec:sphere}

In this section we apply the GFCN machinery to prove Theorem~\ref{thm_sphere}, which states that $m_{\perp}(\Sphere^2) \leq 7/11$.

The problem of avoiding orthonormal bases is closely tied to the Kochen--Specker theorem.
The theorem states that there is no assignment of values $0$ and $1$ to all points of $\Sphere^2$ such that antipodal points receive the same value and every orthonormal basis of $\R^3$ has exactly one vector assigned $1$~\cite{KochenSpecker1967}.
The condition on orthonormal bases is equivalent to requiring that no two orthogonal points both receive $1$ and that no orthonormal basis consists entirely of points receiving $0$.
Our problem retains only the latter restriction, requiring the measurable set under consideration to contain no orthonormal basis without requiring its complement to avoid orthogonal pairs.

Let $\sigma$ denote normalized spherical measure.
If $A\subseteq\Sphere^2$ contains no orthonormal basis, then at most two vectors of any orthonormal basis lie in $A$.
Averaging over uniformly random orthonormal bases, each vector lies in $A$ with probability $\sigma(A)$, so $3\sigma(A)\leq 2$.
Therefore
\begin{equation}\label{eq:trivial-averaging-sphere}
m_\perp(\Sphere^2)\leq\frac{2}{3}.
\end{equation}

For lower bounds, it is useful to consider the classical problem of avoiding orthogonal pairs.
Witsenhausen asked for the largest measurable subset of
\(S^{n-1}\) containing no orthogonal pair, and the expected extremal example is
the double cap: two antipodal open spherical caps of angular radius \(\pi/4\)~\cite{Witsenhausen1974, Kalai}.
On \(\Sphere^2\) a double cap has measure
$1-1/\sqrt 2 \approx 0.2929$. The best known upper bound for Witsenhausen's problem is $0.2953$, obtained in~\cite{CzifraEtAl2026}. Since a single double cap contains no orthogonal pair, the
union of two double caps with orthogonal axes contains no orthogonal triple;
this already gives
\[
    m_\perp(\Sphere^2) \geq 2\left(1-\frac1{\sqrt2}\right)
    =
    2-\sqrt2
    \approx 0.5858.
\]

Recently, Williams and Constantin~\cite{WilliamsConstantin2025} improved this lower bound by constructing a curved
``scalloped belt'' \(N_0\subset S^2\) with no three mutually orthogonal
directions and measure
\[
    \sigma(N_0)
    =
    1-\frac1{\sqrt2}+\frac{\sqrt2\log 2}{\pi}
    \approx
    0.6049.
\]
Their construction supersedes earlier improvements by Perez~\cite{Perez2019}.
They conjecture it to be optimal.
This construction arose from the study of subsets of $\Sphere^2$ admitting assignments of $0$ and $1$ that agree at antipodal points, assign $1$ to no orthogonal pair, and assign exactly one $1$ to every orthonormal basis contained in the subset.
The set $N_0$ is disjoint from the double cap centered at the north and south poles, and their union admits such an assignment, with value $0$ on $N_0$ and value $1$ on the double cap.

Since changing the signs of the vectors in an orthonormal basis independently preserves orthonormality, the set $A\cup(-A)$ is basis-avoiding whenever $A\subseteq\Sphere^2$ is, and $\sigma(A\cup(-A))\geq\sigma(A)$.
We may therefore restrict to antipodally symmetric sets and pass to the projective plane by identifying antipodal points on the sphere.
Let $\pi:\Sphere^2\to\RP^2$ be the antipodal quotient, and endow
$\RP^2$ with the quotient round metric and the probability measure
$\bar\sigma\coloneq\pi_*\sigma$.
Let $\Pcal_\perp$ be the family of projective orthonormal triples
\[
\{\pi(u_1),\pi(u_2),\pi(u_3)\},
\]
where $(u_1,u_2,u_3)$ ranges over orthonormal bases of $\R^3$.
For a probe $X\subseteq\RP^2$, write
$\Gamma_\perp(X)\coloneq\Gamma_{\Pcal_\perp}(X)$.

\begin{lemma}\label{lem:projective-reduction}
We have
\[
m_\perp(\Sphere^2)=m_{\Pcal_\perp}(\RP^2).
\]
\end{lemma}

\begin{proof}
Let $A\subseteq\Sphere^2$ be measurable and contain no orthonormal
basis.
For every $\varepsilon>0$, choose a compact set $K\subseteq A$
with $\sigma(K)\geq\sigma(A)-\varepsilon$.
The set $B\coloneq\pi(K)$ is compact and $\Pcal_\perp$-independent,
since a projective orthonormal triple in $B$ would lift to an
orthonormal basis contained in $K$.
As $\pi^{-1}(B)=K\cup(-K)$, we have
\[
\bar\sigma(B)=\sigma(K\cup(-K))\geq\sigma(A)-\varepsilon.
\]
Taking the supremum over $A$ and letting $\varepsilon\to0$ gives
$m_\perp(\Sphere^2)\leq m_{\Pcal_\perp}(\RP^2)$.

Conversely, if $B\subseteq\RP^2$ is Borel and
$\Pcal_\perp$-independent, then $\pi^{-1}(B)$ is Borel, contains no
orthonormal basis, and satisfies
$\sigma(\pi^{-1}(B))=\bar\sigma(B)$.
This gives the reverse inequality.
\end{proof}

The space $\RP^2$ is compact, and the natural action of $O(3)$ is
continuous, transitive, and isometric.
Both $\bar\sigma$ and $\Pcal_\perp$ are $O(3)$-invariant.
Thus Theorem~\ref{thm:compact-probe-bound} applies to projective probes.
Their random placement is obtained by applying a Haar-uniform element of
$O(3)$ to all the projective points.
Here two projective subconfigurations are congruent precisely when a single
element of $O(3)$ maps one onto the other.

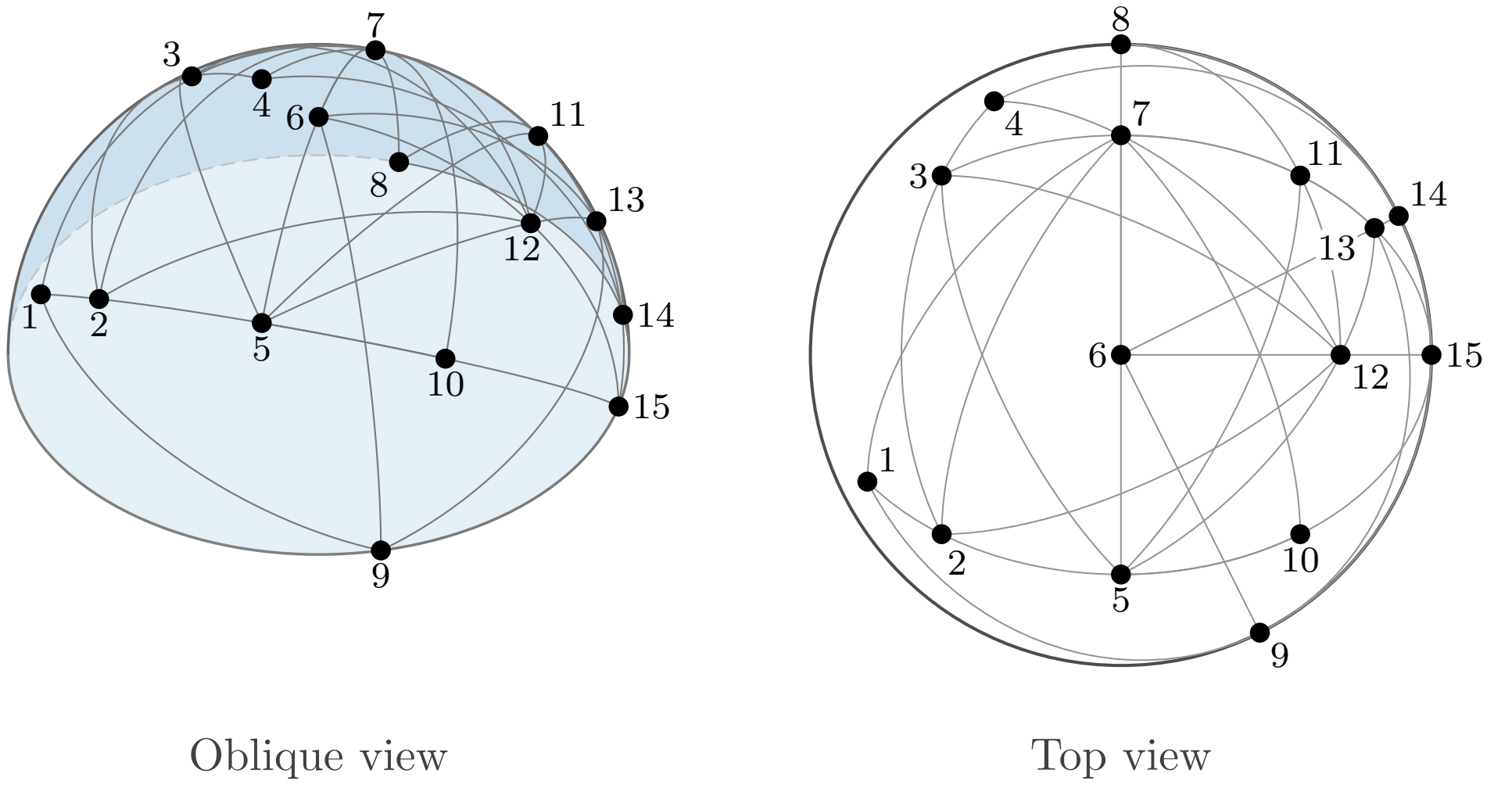
\begin{figure}
\centering
\resizebox{\textwidth}{!}{%
% Two-panel TikZ rendering of the fifteen-point spherical probe.
% For an orthogonal pair u_i,u_j, the arc is cos(t)u_i + sin(t)u_j, 0 <= t <= pi/2.
\begin{tikzpicture}[
    vertex/.style={circle, fill=black, inner sep=1.55pt},
    ortharc/.style={draw=black!42, line width=0.35pt},
    ortharcoblique/.style={draw=black!52, line width=0.37pt},
    numlabel/.style={font=\scriptsize,fill=white,inner sep=0.4pt},
  ]
  % Left panel: orthographic oblique view of the upper hemisphere.
  % Camera azimuth 15 degrees, elevation 40 degrees.
  \begin{scope}[x=2.40cm,y=2.40cm]
    \tikzset{numlabel/.style={font=\scriptsize,inner sep=0.4pt}}
    % Reproduce the original translucent hemisphere shading: the upper crescent
    % receives a second layer of the same tint because the two projected sheets overlap.
    \definecolor{hemisphereBlue}{RGB}{31,119,180}
    \path[fill=hemisphereBlue,fill opacity=0.12] (-1,0) arc[start angle=180,end angle=0,radius=1]
      -- plot[domain=0:180,variable=\t,samples=70]
       ({cos(\t)},{-sin(40)*sin(\t)}) -- cycle;
    \path[fill=hemisphereBlue,fill opacity=0.12] (-1,0) arc[start angle=180,end angle=0,radius=1]
      -- plot[domain=0:180,variable=\t,samples=70]
       ({cos(\t)},{sin(40)*sin(\t)}) -- cycle;
    % Stop the hidden rear equator short of the rim, and place it beneath the outlines.
    \draw[black!24,densely dashed,line width=0.4pt,line cap=round] plot[domain=10:170,variable=\t,samples=70]
       ({cos(\t)},{sin(40)*sin(\t)});
    \draw[black!60,line width=0.6pt] (-1,0) arc[start angle=180,end angle=0,radius=1];
    \draw[black!50,line width=0.6pt] plot[domain=0:180,variable=\t,samples=70]
       ({cos(\t)},{-sin(40)*sin(\t)});
    \coordinate (h1) at ({(-2*cos(15)+-1*sin(15))/sqrt(6)}, {(-sin(40)*sin(15)*(-2)+sin(40)*cos(15)*(-1)+cos(40)*(1))/sqrt(6)});
    \coordinate (h2) at ({(-1*cos(15)+-1*sin(15))/sqrt(3)}, {(-sin(40)*sin(15)*(-1)+sin(40)*cos(15)*(-1)+cos(40)*(1))/sqrt(3)});
    \coordinate (h3) at ({(-1*cos(15)+1*sin(15))/sqrt(3)}, {(-sin(40)*sin(15)*(-1)+sin(40)*cos(15)*(1)+cos(40)*(1))/sqrt(3)});
    \coordinate (h4) at ({(-1*cos(15)+2*sin(15))/sqrt(6)}, {(-sin(40)*sin(15)*(-1)+sin(40)*cos(15)*(2)+cos(40)*(1))/sqrt(6)});
    \coordinate (h5) at ({(0*cos(15)+-1*sin(15))/sqrt(2)}, {(-sin(40)*sin(15)*(0)+sin(40)*cos(15)*(-1)+cos(40)*(1))/sqrt(2)});
    \coordinate (h6) at ({(0*cos(15)+0*sin(15))/sqrt(1)}, {(-sin(40)*sin(15)*(0)+sin(40)*cos(15)*(0)+cos(40)*(1))/sqrt(1)});
    \coordinate (h7) at ({(0*cos(15)+1*sin(15))/sqrt(2)}, {(-sin(40)*sin(15)*(0)+sin(40)*cos(15)*(1)+cos(40)*(1))/sqrt(2)});
    \coordinate (h8) at ({(0*cos(15)+1*sin(15))/sqrt(1)}, {(-sin(40)*sin(15)*(0)+sin(40)*cos(15)*(1)+cos(40)*(0))/sqrt(1)});
    \coordinate (h9) at ({(1*cos(15)+-2*sin(15))/sqrt(5)}, {(-sin(40)*sin(15)*(1)+sin(40)*cos(15)*(-2)+cos(40)*(0))/sqrt(5)});
    \coordinate (h10) at ({(1*cos(15)+-1*sin(15))/sqrt(3)}, {(-sin(40)*sin(15)*(1)+sin(40)*cos(15)*(-1)+cos(40)*(1))/sqrt(3)});
    \coordinate (h11) at ({(1*cos(15)+1*sin(15))/sqrt(3)}, {(-sin(40)*sin(15)*(1)+sin(40)*cos(15)*(1)+cos(40)*(1))/sqrt(3)});
    \coordinate (h12) at ({(1*cos(15)+0*sin(15))/sqrt(2)}, {(-sin(40)*sin(15)*(1)+sin(40)*cos(15)*(0)+cos(40)*(1))/sqrt(2)});
    \coordinate (h13) at ({(2*cos(15)+1*sin(15))/sqrt(6)}, {(-sin(40)*sin(15)*(2)+sin(40)*cos(15)*(1)+cos(40)*(1))/sqrt(6)});
    \coordinate (h14) at ({(2*cos(15)+1*sin(15))/sqrt(5)}, {(-sin(40)*sin(15)*(2)+sin(40)*cos(15)*(1)+cos(40)*(0))/sqrt(5)});
    \coordinate (h15) at ({(1*cos(15)+0*sin(15))/sqrt(1)}, {(-sin(40)*sin(15)*(1)+sin(40)*cos(15)*(0)+cos(40)*(0))/sqrt(1)});
    % Exact oblique projections of the same 25 quarter great-circle arcs.
    \draw[ortharcoblique] plot[domain=0:90,variable=\t,samples=36]
      ({(-2*cos(15)+-1*sin(15))*cos(\t)/sqrt(6)+(0*cos(15)+1*sin(15))*sin(\t)/sqrt(2)},
       {(-sin(40)*sin(15)*(-2)+sin(40)*cos(15)*(-1)+cos(40)*(1))*cos(\t)/sqrt(6)+(-sin(40)*sin(15)*(0)+sin(40)*cos(15)*(1)+cos(40)*(1))*sin(\t)/sqrt(2)});
    \draw[ortharcoblique] plot[domain=0:90,variable=\t,samples=36]
      ({(-2*cos(15)+-1*sin(15))*cos(\t)/sqrt(6)+(1*cos(15)+-2*sin(15))*sin(\t)/sqrt(5)},
       {(-sin(40)*sin(15)*(-2)+sin(40)*cos(15)*(-1)+cos(40)*(1))*cos(\t)/sqrt(6)+(-sin(40)*sin(15)*(1)+sin(40)*cos(15)*(-2)+cos(40)*(0))*sin(\t)/sqrt(5)});
    \draw[ortharcoblique] plot[domain=0:90,variable=\t,samples=36]
      ({(-2*cos(15)+-1*sin(15))*cos(\t)/sqrt(6)+(1*cos(15)+-1*sin(15))*sin(\t)/sqrt(3)},
       {(-sin(40)*sin(15)*(-2)+sin(40)*cos(15)*(-1)+cos(40)*(1))*cos(\t)/sqrt(6)+(-sin(40)*sin(15)*(1)+sin(40)*cos(15)*(-1)+cos(40)*(1))*sin(\t)/sqrt(3)});
    \draw[ortharcoblique] plot[domain=0:90,variable=\t,samples=36]
      ({(-1*cos(15)+-1*sin(15))*cos(\t)/sqrt(3)+(-1*cos(15)+2*sin(15))*sin(\t)/sqrt(6)},
       {(-sin(40)*sin(15)*(-1)+sin(40)*cos(15)*(-1)+cos(40)*(1))*cos(\t)/sqrt(3)+(-sin(40)*sin(15)*(-1)+sin(40)*cos(15)*(2)+cos(40)*(1))*sin(\t)/sqrt(6)});
    \draw[ortharcoblique] plot[domain=0:90,variable=\t,samples=36]
      ({(-1*cos(15)+-1*sin(15))*cos(\t)/sqrt(3)+(0*cos(15)+1*sin(15))*sin(\t)/sqrt(2)},
       {(-sin(40)*sin(15)*(-1)+sin(40)*cos(15)*(-1)+cos(40)*(1))*cos(\t)/sqrt(3)+(-sin(40)*sin(15)*(0)+sin(40)*cos(15)*(1)+cos(40)*(1))*sin(\t)/sqrt(2)});
    \draw[ortharcoblique] plot[domain=0:90,variable=\t,samples=36]
      ({(-1*cos(15)+-1*sin(15))*cos(\t)/sqrt(3)+(1*cos(15)+0*sin(15))*sin(\t)/sqrt(2)},
       {(-sin(40)*sin(15)*(-1)+sin(40)*cos(15)*(-1)+cos(40)*(1))*cos(\t)/sqrt(3)+(-sin(40)*sin(15)*(1)+sin(40)*cos(15)*(0)+cos(40)*(1))*sin(\t)/sqrt(2)});
    \draw[ortharcoblique] plot[domain=0:90,variable=\t,samples=36]
      ({(-1*cos(15)+1*sin(15))*cos(\t)/sqrt(3)+(0*cos(15)+-1*sin(15))*sin(\t)/sqrt(2)},
       {(-sin(40)*sin(15)*(-1)+sin(40)*cos(15)*(1)+cos(40)*(1))*cos(\t)/sqrt(3)+(-sin(40)*sin(15)*(0)+sin(40)*cos(15)*(-1)+cos(40)*(1))*sin(\t)/sqrt(2)});
    \draw[ortharcoblique] plot[domain=0:90,variable=\t,samples=36]
      ({(-1*cos(15)+1*sin(15))*cos(\t)/sqrt(3)+(1*cos(15)+0*sin(15))*sin(\t)/sqrt(2)},
       {(-sin(40)*sin(15)*(-1)+sin(40)*cos(15)*(1)+cos(40)*(1))*cos(\t)/sqrt(3)+(-sin(40)*sin(15)*(1)+sin(40)*cos(15)*(0)+cos(40)*(1))*sin(\t)/sqrt(2)});
    \draw[ortharcoblique] plot[domain=0:90,variable=\t,samples=36]
      ({(-1*cos(15)+1*sin(15))*cos(\t)/sqrt(3)+(2*cos(15)+1*sin(15))*sin(\t)/sqrt(6)},
       {(-sin(40)*sin(15)*(-1)+sin(40)*cos(15)*(1)+cos(40)*(1))*cos(\t)/sqrt(3)+(-sin(40)*sin(15)*(2)+sin(40)*cos(15)*(1)+cos(40)*(1))*sin(\t)/sqrt(6)});
    \draw[ortharcoblique] plot[domain=0:90,variable=\t,samples=36]
      ({(-1*cos(15)+2*sin(15))*cos(\t)/sqrt(6)+(1*cos(15)+0*sin(15))*sin(\t)/sqrt(2)},
       {(-sin(40)*sin(15)*(-1)+sin(40)*cos(15)*(2)+cos(40)*(1))*cos(\t)/sqrt(6)+(-sin(40)*sin(15)*(1)+sin(40)*cos(15)*(0)+cos(40)*(1))*sin(\t)/sqrt(2)});
    \draw[ortharcoblique] plot[domain=0:90,variable=\t,samples=36]
      ({(-1*cos(15)+2*sin(15))*cos(\t)/sqrt(6)+(2*cos(15)+1*sin(15))*sin(\t)/sqrt(5)},
       {(-sin(40)*sin(15)*(-1)+sin(40)*cos(15)*(2)+cos(40)*(1))*cos(\t)/sqrt(6)+(-sin(40)*sin(15)*(2)+sin(40)*cos(15)*(1)+cos(40)*(0))*sin(\t)/sqrt(5)});
    \draw[ortharcoblique] plot[domain=0:90,variable=\t,samples=36]
      ({(0*cos(15)+-1*sin(15))*cos(\t)/sqrt(2)+(0*cos(15)+1*sin(15))*sin(\t)/sqrt(2)},
       {(-sin(40)*sin(15)*(0)+sin(40)*cos(15)*(-1)+cos(40)*(1))*cos(\t)/sqrt(2)+(-sin(40)*sin(15)*(0)+sin(40)*cos(15)*(1)+cos(40)*(1))*sin(\t)/sqrt(2)});
    \draw[ortharcoblique] plot[domain=0:90,variable=\t,samples=36]
      ({(0*cos(15)+-1*sin(15))*cos(\t)/sqrt(2)+(1*cos(15)+1*sin(15))*sin(\t)/sqrt(3)},
       {(-sin(40)*sin(15)*(0)+sin(40)*cos(15)*(-1)+cos(40)*(1))*cos(\t)/sqrt(2)+(-sin(40)*sin(15)*(1)+sin(40)*cos(15)*(1)+cos(40)*(1))*sin(\t)/sqrt(3)});
    \draw[ortharcoblique] plot[domain=0:90,variable=\t,samples=36]
      ({(0*cos(15)+-1*sin(15))*cos(\t)/sqrt(2)+(2*cos(15)+1*sin(15))*sin(\t)/sqrt(6)},
       {(-sin(40)*sin(15)*(0)+sin(40)*cos(15)*(-1)+cos(40)*(1))*cos(\t)/sqrt(2)+(-sin(40)*sin(15)*(2)+sin(40)*cos(15)*(1)+cos(40)*(1))*sin(\t)/sqrt(6)});
    \draw[ortharcoblique] plot[domain=0:90,variable=\t,samples=36]
      ({(0*cos(15)+-1*sin(15))*cos(\t)/sqrt(2)+(1*cos(15)+0*sin(15))*sin(\t)/sqrt(1)},
       {(-sin(40)*sin(15)*(0)+sin(40)*cos(15)*(-1)+cos(40)*(1))*cos(\t)/sqrt(2)+(-sin(40)*sin(15)*(1)+sin(40)*cos(15)*(0)+cos(40)*(0))*sin(\t)/sqrt(1)});
    \draw[ortharcoblique] plot[domain=0:90,variable=\t,samples=36]
      ({(0*cos(15)+0*sin(15))*cos(\t)/sqrt(1)+(0*cos(15)+1*sin(15))*sin(\t)/sqrt(1)},
       {(-sin(40)*sin(15)*(0)+sin(40)*cos(15)*(0)+cos(40)*(1))*cos(\t)/sqrt(1)+(-sin(40)*sin(15)*(0)+sin(40)*cos(15)*(1)+cos(40)*(0))*sin(\t)/sqrt(1)});
    \draw[ortharcoblique] plot[domain=0:90,variable=\t,samples=36]
      ({(0*cos(15)+0*sin(15))*cos(\t)/sqrt(1)+(1*cos(15)+-2*sin(15))*sin(\t)/sqrt(5)},
       {(-sin(40)*sin(15)*(0)+sin(40)*cos(15)*(0)+cos(40)*(1))*cos(\t)/sqrt(1)+(-sin(40)*sin(15)*(1)+sin(40)*cos(15)*(-2)+cos(40)*(0))*sin(\t)/sqrt(5)});
    \draw[ortharcoblique] plot[domain=0:90,variable=\t,samples=36]
      ({(0*cos(15)+0*sin(15))*cos(\t)/sqrt(1)+(2*cos(15)+1*sin(15))*sin(\t)/sqrt(5)},
       {(-sin(40)*sin(15)*(0)+sin(40)*cos(15)*(0)+cos(40)*(1))*cos(\t)/sqrt(1)+(-sin(40)*sin(15)*(2)+sin(40)*cos(15)*(1)+cos(40)*(0))*sin(\t)/sqrt(5)});
    \draw[ortharcoblique] plot[domain=0:90,variable=\t,samples=36]
      ({(0*cos(15)+0*sin(15))*cos(\t)/sqrt(1)+(1*cos(15)+0*sin(15))*sin(\t)/sqrt(1)},
       {(-sin(40)*sin(15)*(0)+sin(40)*cos(15)*(0)+cos(40)*(1))*cos(\t)/sqrt(1)+(-sin(40)*sin(15)*(1)+sin(40)*cos(15)*(0)+cos(40)*(0))*sin(\t)/sqrt(1)});
    \draw[ortharcoblique] plot[domain=0:90,variable=\t,samples=36]
      ({(0*cos(15)+1*sin(15))*cos(\t)/sqrt(2)+(1*cos(15)+-1*sin(15))*sin(\t)/sqrt(3)},
       {(-sin(40)*sin(15)*(0)+sin(40)*cos(15)*(1)+cos(40)*(1))*cos(\t)/sqrt(2)+(-sin(40)*sin(15)*(1)+sin(40)*cos(15)*(-1)+cos(40)*(1))*sin(\t)/sqrt(3)});
    \draw[ortharcoblique] plot[domain=0:90,variable=\t,samples=36]
      ({(0*cos(15)+1*sin(15))*cos(\t)/sqrt(2)+(1*cos(15)+0*sin(15))*sin(\t)/sqrt(1)},
       {(-sin(40)*sin(15)*(0)+sin(40)*cos(15)*(1)+cos(40)*(1))*cos(\t)/sqrt(2)+(-sin(40)*sin(15)*(1)+sin(40)*cos(15)*(0)+cos(40)*(0))*sin(\t)/sqrt(1)});
    \draw[ortharcoblique] plot[domain=0:90,variable=\t,samples=36]
      ({(0*cos(15)+1*sin(15))*cos(\t)/sqrt(1)+(1*cos(15)+0*sin(15))*sin(\t)/sqrt(2)},
       {(-sin(40)*sin(15)*(0)+sin(40)*cos(15)*(1)+cos(40)*(0))*cos(\t)/sqrt(1)+(-sin(40)*sin(15)*(1)+sin(40)*cos(15)*(0)+cos(40)*(1))*sin(\t)/sqrt(2)});
    \draw[ortharcoblique] plot[domain=0:90,variable=\t,samples=36]
      ({(0*cos(15)+1*sin(15))*cos(\t)/sqrt(1)+(1*cos(15)+0*sin(15))*sin(\t)/sqrt(1)},
       {(-sin(40)*sin(15)*(0)+sin(40)*cos(15)*(1)+cos(40)*(0))*cos(\t)/sqrt(1)+(-sin(40)*sin(15)*(1)+sin(40)*cos(15)*(0)+cos(40)*(0))*sin(\t)/sqrt(1)});
    \draw[ortharcoblique] plot[domain=0:90,variable=\t,samples=36]
      ({(1*cos(15)+-2*sin(15))*cos(\t)/sqrt(5)+(2*cos(15)+1*sin(15))*sin(\t)/sqrt(6)},
       {(-sin(40)*sin(15)*(1)+sin(40)*cos(15)*(-2)+cos(40)*(0))*cos(\t)/sqrt(5)+(-sin(40)*sin(15)*(2)+sin(40)*cos(15)*(1)+cos(40)*(1))*sin(\t)/sqrt(6)});
    \draw[ortharcoblique] plot[domain=0:90,variable=\t,samples=36]
      ({(1*cos(15)+-2*sin(15))*cos(\t)/sqrt(5)+(2*cos(15)+1*sin(15))*sin(\t)/sqrt(5)},
       {(-sin(40)*sin(15)*(1)+sin(40)*cos(15)*(-2)+cos(40)*(0))*cos(\t)/sqrt(5)+(-sin(40)*sin(15)*(2)+sin(40)*cos(15)*(1)+cos(40)*(0))*sin(\t)/sqrt(5)});
    \node[vertex,label={[numlabel,xshift=2pt]below left:$1$}] at (h1) {};
    \node[vertex,label={[numlabel]below:$2$}] at (h2) {};
    \node[vertex,label={[numlabel]above left:$3$}] at (h3) {};
    \node[vertex,label={[numlabel]below:$4$}] at (h4) {};
    \node[vertex,label={[numlabel]below:$5$}] at (h5) {};
    \node[vertex,label={[numlabel]left:$6$}] at (h6) {};
    \node[vertex,label={[numlabel]above:$7$}] at (h7) {};
    \node[vertex,label={[numlabel]below left:$8$}] at (h8) {};
    \node[vertex,label={[numlabel]below:$9$}] at (h9) {};
    \node[vertex,label={[numlabel]below:$10$}] at (h10) {};
    \node[vertex,label={[numlabel]above right:$11$}] at (h11) {};
    \node[vertex,label={[numlabel,xshift=-2pt]below:$12$}] at (h12) {};
    \node[vertex,label={[numlabel]above right:$13$}] at (h13) {};
    \node[vertex,label={[numlabel]right:$14$}] at (h14) {};
    \node[vertex,label={[numlabel]right:$15$}] at (h15) {};
    \node[font=\small, text=black!75] at (0,-1.30) {Oblique view};
  \end{scope}

  % Right panel: the xy-projection of the exact quarter great-circle arcs.
  % Match the sphere diameter in the two views.
  \begin{scope}[shift={(6.2cm,0)},x=2.40cm,y=2.40cm]
    \tikzset{
      numlabel/.style={font=\scriptsize,inner sep=0.4pt},
      cutlabel/.style={font=\scriptsize,shape=circle,fill=white,inner sep=0.1pt,outer sep=0pt}
    }
    \draw[black!70,line width=0.7pt] (0,0) circle (1);
    \coordinate (p1) at ({(-2)/sqrt(6)},{(-1)/sqrt(6)});
    \coordinate (p2) at ({(-1)/sqrt(3)},{(-1)/sqrt(3)});
    \coordinate (p3) at ({(-1)/sqrt(3)},{(1)/sqrt(3)});
    \coordinate (p4) at ({(-1)/sqrt(6)},{(2)/sqrt(6)});
    \coordinate (p5) at ({(0)/sqrt(2)},{(-1)/sqrt(2)});
    \coordinate (p6) at ({(0)/sqrt(1)},{(0)/sqrt(1)});
    \coordinate (p7) at ({(0)/sqrt(2)},{(1)/sqrt(2)});
    \coordinate (p8) at ({(0)/sqrt(1)},{(1)/sqrt(1)});
    \coordinate (p9) at ({(1)/sqrt(5)},{(-2)/sqrt(5)});
    \coordinate (p10) at ({(1)/sqrt(3)},{(-1)/sqrt(3)});
    \coordinate (p11) at ({(1)/sqrt(3)},{(1)/sqrt(3)});
    \coordinate (p12) at ({(1)/sqrt(2)},{(0)/sqrt(2)});
    \coordinate (p13) at ({(2)/sqrt(6)},{(1)/sqrt(6)});
    \coordinate (p14) at ({(2)/sqrt(5)},{(1)/sqrt(5)});
    \coordinate (p15) at ({(1)/sqrt(1)},{(0)/sqrt(1)});
    % Exact projected great-circle quarter arcs, one for each orthogonal pair.
    \draw[ortharc] plot[domain=0:90,variable=\t,samples=28]
      ({(-2)*cos(\t)/sqrt(6)+(0)*sin(\t)/sqrt(2)},
       {(-1)*cos(\t)/sqrt(6)+(1)*sin(\t)/sqrt(2)});
    \draw[ortharc] plot[domain=0:90,variable=\t,samples=28]
      ({(-2)*cos(\t)/sqrt(6)+(1)*sin(\t)/sqrt(5)},
       {(-1)*cos(\t)/sqrt(6)+(-2)*sin(\t)/sqrt(5)});
    \draw[ortharc] plot[domain=0:90,variable=\t,samples=28]
      ({(-2)*cos(\t)/sqrt(6)+(1)*sin(\t)/sqrt(3)},
       {(-1)*cos(\t)/sqrt(6)+(-1)*sin(\t)/sqrt(3)});
    \draw[ortharc] plot[domain=0:90,variable=\t,samples=28]
      ({(-1)*cos(\t)/sqrt(3)+(-1)*sin(\t)/sqrt(6)},
       {(-1)*cos(\t)/sqrt(3)+(2)*sin(\t)/sqrt(6)});
    \draw[ortharc] plot[domain=0:90,variable=\t,samples=28]
      ({(-1)*cos(\t)/sqrt(3)+(0)*sin(\t)/sqrt(2)},
       {(-1)*cos(\t)/sqrt(3)+(1)*sin(\t)/sqrt(2)});
    \draw[ortharc] plot[domain=0:90,variable=\t,samples=28]
      ({(-1)*cos(\t)/sqrt(3)+(1)*sin(\t)/sqrt(2)},
       {(-1)*cos(\t)/sqrt(3)+(0)*sin(\t)/sqrt(2)});
    \draw[ortharc] plot[domain=0:90,variable=\t,samples=28]
      ({(-1)*cos(\t)/sqrt(3)+(0)*sin(\t)/sqrt(2)},
       {(1)*cos(\t)/sqrt(3)+(-1)*sin(\t)/sqrt(2)});
    \draw[ortharc] plot[domain=0:90,variable=\t,samples=28]
      ({(-1)*cos(\t)/sqrt(3)+(1)*sin(\t)/sqrt(2)},
       {(1)*cos(\t)/sqrt(3)+(0)*sin(\t)/sqrt(2)});
    \draw[ortharc] plot[domain=0:90,variable=\t,samples=28]
      ({(-1)*cos(\t)/sqrt(3)+(2)*sin(\t)/sqrt(6)},
       {(1)*cos(\t)/sqrt(3)+(1)*sin(\t)/sqrt(6)});
    \draw[ortharc] plot[domain=0:90,variable=\t,samples=28]
      ({(-1)*cos(\t)/sqrt(6)+(1)*sin(\t)/sqrt(2)},
       {(2)*cos(\t)/sqrt(6)+(0)*sin(\t)/sqrt(2)});
    \draw[ortharc] plot[domain=0:90,variable=\t,samples=28]
      ({(-1)*cos(\t)/sqrt(6)+(2)*sin(\t)/sqrt(5)},
       {(2)*cos(\t)/sqrt(6)+(1)*sin(\t)/sqrt(5)});
    \draw[ortharc] plot[domain=0:90,variable=\t,samples=28]
      ({(0)*cos(\t)/sqrt(2)+(0)*sin(\t)/sqrt(2)},
       {(-1)*cos(\t)/sqrt(2)+(1)*sin(\t)/sqrt(2)});
    \draw[ortharc] plot[domain=0:90,variable=\t,samples=28]
      ({(0)*cos(\t)/sqrt(2)+(1)*sin(\t)/sqrt(3)},
       {(-1)*cos(\t)/sqrt(2)+(1)*sin(\t)/sqrt(3)});
    \draw[ortharc] plot[domain=0:90,variable=\t,samples=28]
      ({(0)*cos(\t)/sqrt(2)+(2)*sin(\t)/sqrt(6)},
       {(-1)*cos(\t)/sqrt(2)+(1)*sin(\t)/sqrt(6)});
    \draw[ortharc] plot[domain=0:90,variable=\t,samples=28]
      ({(0)*cos(\t)/sqrt(2)+(1)*sin(\t)/sqrt(1)},
       {(-1)*cos(\t)/sqrt(2)+(0)*sin(\t)/sqrt(1)});
    \draw[ortharc] plot[domain=0:90,variable=\t,samples=28]
      ({(0)*cos(\t)/sqrt(1)+(0)*sin(\t)/sqrt(1)},
       {(0)*cos(\t)/sqrt(1)+(1)*sin(\t)/sqrt(1)});
    \draw[ortharc] plot[domain=0:90,variable=\t,samples=28]
      ({(0)*cos(\t)/sqrt(1)+(1)*sin(\t)/sqrt(5)},
       {(0)*cos(\t)/sqrt(1)+(-2)*sin(\t)/sqrt(5)});
    \draw[ortharc] plot[domain=0:90,variable=\t,samples=28]
      ({(0)*cos(\t)/sqrt(1)+(2)*sin(\t)/sqrt(5)},
       {(0)*cos(\t)/sqrt(1)+(1)*sin(\t)/sqrt(5)});
    \draw[ortharc] plot[domain=0:90,variable=\t,samples=28]
      ({(0)*cos(\t)/sqrt(1)+(1)*sin(\t)/sqrt(1)},
       {(0)*cos(\t)/sqrt(1)+(0)*sin(\t)/sqrt(1)});
    \draw[ortharc] plot[domain=0:90,variable=\t,samples=28]
      ({(0)*cos(\t)/sqrt(2)+(1)*sin(\t)/sqrt(3)},
       {(1)*cos(\t)/sqrt(2)+(-1)*sin(\t)/sqrt(3)});
    \draw[ortharc] plot[domain=0:90,variable=\t,samples=28]
      ({(0)*cos(\t)/sqrt(2)+(1)*sin(\t)/sqrt(1)},
       {(1)*cos(\t)/sqrt(2)+(0)*sin(\t)/sqrt(1)});
    \draw[ortharc] plot[domain=0:90,variable=\t,samples=28]
      ({(0)*cos(\t)/sqrt(1)+(1)*sin(\t)/sqrt(2)},
       {(1)*cos(\t)/sqrt(1)+(0)*sin(\t)/sqrt(2)});
    \draw[ortharc] plot[domain=0:90,variable=\t,samples=28]
      ({(0)*cos(\t)/sqrt(1)+(1)*sin(\t)/sqrt(1)},
       {(1)*cos(\t)/sqrt(1)+(0)*sin(\t)/sqrt(1)});
    \draw[ortharc] plot[domain=0:90,variable=\t,samples=28]
      ({(1)*cos(\t)/sqrt(5)+(2)*sin(\t)/sqrt(6)},
       {(-2)*cos(\t)/sqrt(5)+(1)*sin(\t)/sqrt(6)});
    \draw[ortharc] plot[domain=0:90,variable=\t,samples=28]
      ({(1)*cos(\t)/sqrt(5)+(2)*sin(\t)/sqrt(5)},
       {(-2)*cos(\t)/sqrt(5)+(1)*sin(\t)/sqrt(5)});
    \node[vertex,label={[numlabel]above right:$1$}] at (p1) {};
    \node[vertex,label={[numlabel,xshift=-1,yshift=-1.5pt]below right:$2$}] at (p2) {};
    \node[vertex,label={[numlabel]left:$3$}] at (p3) {};
    \node[vertex,label={[numlabel]below right:$4$}] at (p4) {};
    \node[vertex,label={[numlabel]below:$5$}] at (p5) {};
    \node[vertex,label={[numlabel]left:$6$}] at (p6) {};
    \node[vertex,label={[numlabel]above right:$7$}] at (p7) {};
    \node[vertex,label={[numlabel]above:$8$}] at (p8) {};
    \node[vertex,label={[numlabel]below right:$9$}] at (p9) {};
    \node[vertex,label={[numlabel]below:$10$}] at (p10) {};
    \node[vertex,label={[numlabel,xshift=-1pt]above right:$11$}] at (p11) {};
    \node[vertex,label={[numlabel]below right:$12$}] at (p12) {};
    \node[vertex,label={[cutlabel,xshift=-3pt,yshift=1pt]below left:$13$}] at (p13) {};
    \node[vertex,label={[numlabel]above right:$14$}] at (p14) {};
    \node[vertex,label={[numlabel]right:$15$}] at (p15) {};
    \node[font=\small,text=black!75] at (0,-1.30) {Top view};
  \end{scope}
\end{tikzpicture}
}
\caption{The fifteen-point projective probe $X_{15}$, represented by unit vectors in the closed northern hemisphere of $\Sphere^2$.
Both panels show the same labeled representatives and quarter great-circle arcs joining orthogonal pairs.}
\label{fig:north-hemisphere-points}
\end{figure}

We now present a fifteen-point probe $X_{15}\subseteq\RP^2$.
The fifteen points of $X_{15}$ are represented by the nonzero integer vectors $(x_j,y_j,z_j)$ listed in Table~\ref{tab:vertices-projective}.
Figure~\ref{fig:north-hemisphere-points} shows the corresponding unit representatives chosen in the closed northern hemisphere.
The left panel gives an oblique view of the hemisphere, while the right panel shows its orthogonal projection onto the plane $z=0$.
Quarter great-circle arcs indicate orthogonal pairs among the chosen representatives.
The labels in both panels match the row order of Table~\ref{tab:vertices-projective}.

\begin{table}
\centering
\small
\begin{tabular}{c|rrr}
\hline
$j$ & $x_j$ & $y_j$ & $z_j$ \\
\hline
1 & $-2$ & $-1$ & $1$ \\
2 & $-1$ & $-1$ & $1$ \\
3 & $-1$ & $1$ & $1$ \\
4 & $-1$ & $2$ & $1$ \\
5 & $0$ & $-1$ & $1$ \\
6 & $0$ & $0$ & $1$ \\
7 & $0$ & $1$ & $1$ \\
8 & $0$ & $1$ & $0$ \\
9 & $1$ & $-2$ & $0$ \\
10 & $1$ & $-1$ & $1$ \\
11 & $1$ & $1$ & $1$ \\
12 & $1$ & $0$ & $1$ \\
13 & $2$ & $1$ & $1$ \\
14 & $2$ & $1$ & $0$ \\
15 & $1$ & $0$ & $0$ \\
\end{tabular}
\vspace{0.3cm}
\caption{Integer vectors representing the projective points of $X_{15}$.
The indices match the labels in Figure~\ref{fig:north-hemisphere-points}.}
\label{tab:vertices-projective}
\end{table}

For this probe, the finite linear program gives
$\Gamma_\perp(X_{15})=7/11$.
Combining Lemma~\ref{lem:projective-reduction} and
Theorem~\ref{thm:compact-probe-bound}, we obtain
\[
m_\perp(\Sphere^2)=m_{\Pcal_\perp}(\RP^2)
\leq\Gamma_\perp(X_{15})=\frac{7}{11},
\]
as asserted in Theorem~\ref{thm_sphere}.

\section{Unit equilateral triangles in the plane}\label{sec:plane}

The seminal paper by Erdős et al.~\cite{ERDOS1973341}, which introduced Euclidean Ramsey theory, begins with a simple illustrative example:

\medskip

\begin{quote}
We are concerned in this paper with problems of the type illustrated as follows: Is it true that for any partition of the Euclidean plane into two classes (we say that the plane is \emph{two-colored}), there exists a set of three points all in the same class forming the vertices of an equilateral triangle of side length 1? (We call such a set \emph{monochromatic}.)

In this example the answer is ``no,'' as can be seen by dividing the points $(x, y)$ into two classes according to the parity of $[2y/\sqrt{3}]$.
\end{quote}

% TODO NE ÉGJÜNK ITT BE TRIVIÁLISAN HAMIS VAGY TRIVIÁLISAN IGAZ ÁLLÍTÁSOKKAL.
% Thus, \cite{ERDOS1973341} observes that the plane can be two-colored such that each color class is unit equilateral triangle free ($\Delta$-free, for short). Interestingly, for all known such two-colorings, the upper density of both color classes is $1/2$. Even in the independence density variant, where we allow the second color class to contain $\Delta$, we still apparently cannot exploit this to bring the density of the first color class above $1/2$.

The two color classes in this construction both avoid unit equilateral
triangles and have upper Banach density $1/2$. A different family of ``zebra-stripe'' constructions was given by Jelínek et al. \cite{Jelinek2009}, achieving the same density. In contrast, Currier, Moore, and Yip \cite{CurrierMooreYip2024} showed that every two-coloring of the plane contains a monochromatic congruent copy of any prescribed three-term arithmetic progression.

Let \(\Delta\subseteq\R^2\) be the vertex set of a unit equilateral triangle,
and let
\[
\Pcal_\Delta\coloneq\{R\Delta+v : R\in O(2),\ v\in\R^2\}.
\]
We write
\[
m_\Delta(\R^2)\coloneq m_{\Pcal_\Delta}(\R^2)
\qquad\text{and}\qquad
\Gamma_\Delta(X)\coloneq\Gamma_{\Pcal_\Delta}(X).
\]
Thus $m_\Delta(\R^2)\geq 1/2$.
Although the independence density problem imposes no restriction
on the complement of a $\Pcal_\Delta$-independent set,
no such set with upper Banach density greater than $1/2$
is known to us.
We conjecture that $m_\Delta(\R^2)=1/2$.

To describe the two probes used below, we identify $\R^2$ with the complex plane and set $\omega=(1+i\sqrt{3})/2$ and $\eta=(5+i\sqrt{11})/6$.

We first establish an elementary upper bound using a variation of the construction of Leo and William Moser~\cite{Moser1961}.

\begin{theorem}\label{thm_moser_like}
$m_\Delta(\R^2) \leq 5/8.$
\end{theorem}

\begin{proof}
Adding one vertex to the Moser spindle gives the eight-vertex unit-distance graph $M'$ shown in Figure~\ref{fig_moser_prime}.

\begin{figure}
  \centering
  \begin{tikzpicture}[
    x=3cm,
    y=3cm,
    vertex/.style={circle, fill=black, inner sep=1.8pt}
  ]
    \coordinate (v1) at (0, 0);
    \coordinate (v2) at (1, 0);
    \coordinate (v3) at (0.5, {sqrt(3)/2});
    \coordinate (v4) at (1.5, {sqrt(3)/2});
    \coordinate (v5) at ({5/6}, {sqrt(11)/6});
    \coordinate (v6) at ({(5-sqrt(33))/12}, {(sqrt(11)+5*sqrt(3))/12});
    \coordinate (v7) at ({(15-sqrt(33))/12}, {(3*sqrt(11)+5*sqrt(3))/12});
    \coordinate (v8) at ({(9-sqrt(33))/6}, {sqrt(3)/3});

    \draw[thick]
      (v1) -- (v2) -- (v3) -- (v1)
      (v2) -- (v4) -- (v3)
      (v1) -- (v5) -- (v6) -- (v1)
      (v5) -- (v7) -- (v6)
      (v4) -- (v7) -- (v8) -- (v4);

    \node[vertex, label={[font=\small]below left:$1$}] at (v1) {};
    \node[vertex, label={[font=\small]below right:$2$}] at (v2) {};
    \node[vertex, label={[font=\small]above left:$3$}] at (v3) {};
    \node[vertex, label={[font=\small]right:$4$}] at (v4) {};
    \node[vertex, label={[font=\small]right:$5$}] at (v5) {};
    \node[vertex, label={[font=\small]left:$6$}] at (v6) {};
    \node[vertex, label={[font=\small]above:$7$}] at (v7) {};
    \node[vertex, label={[font=\small]below left:$8$}] at (v8) {};
  \end{tikzpicture}
\caption{The eight-point probe $X_8$ and its unit-distance graph $M'$, obtained by adding vertex $8$ to the Moser spindle.}
\label{fig_moser_prime}
\end{figure}
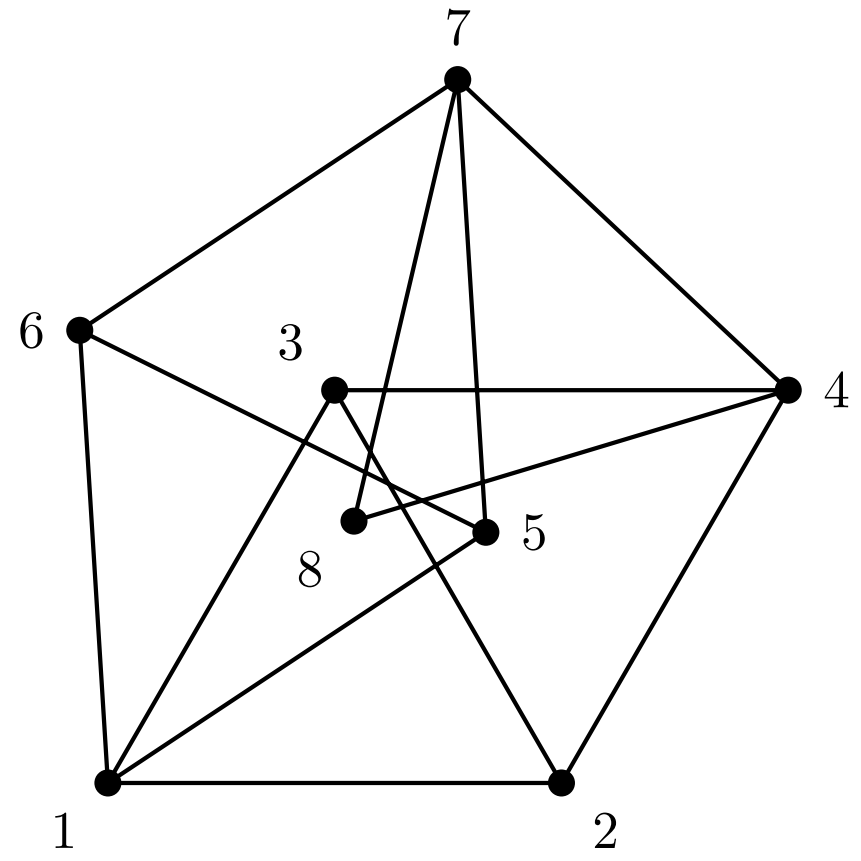

It is easy to verify that any six vertices of $M'$ induce a subgraph containing a triangle. By Lemma~\ref{lem:periodic-reduction-euclidean}, it suffices to prove the bound for $L$-periodic Borel triangle-free sets $A$.
By this property of $M'$, every placement of $X_8$ in the plane has at most five vertices in $A$.
By linearity of expectation, this implies that $\delta_L(A)\leq 5/8$.
\end{proof}

The $5/8$ bound uses only that every triangle-free subset of $X_8$ has at most five points and that each point of $X_8$ has the same inclusion probability under a random placement.
The GFCN linear program also imposes congruence constraints for subconfigurations containing more than one point.
We apply these additional constraints to a fourteen-point probe $X_{14}=\{z_1,\dots,z_{14}\}$.

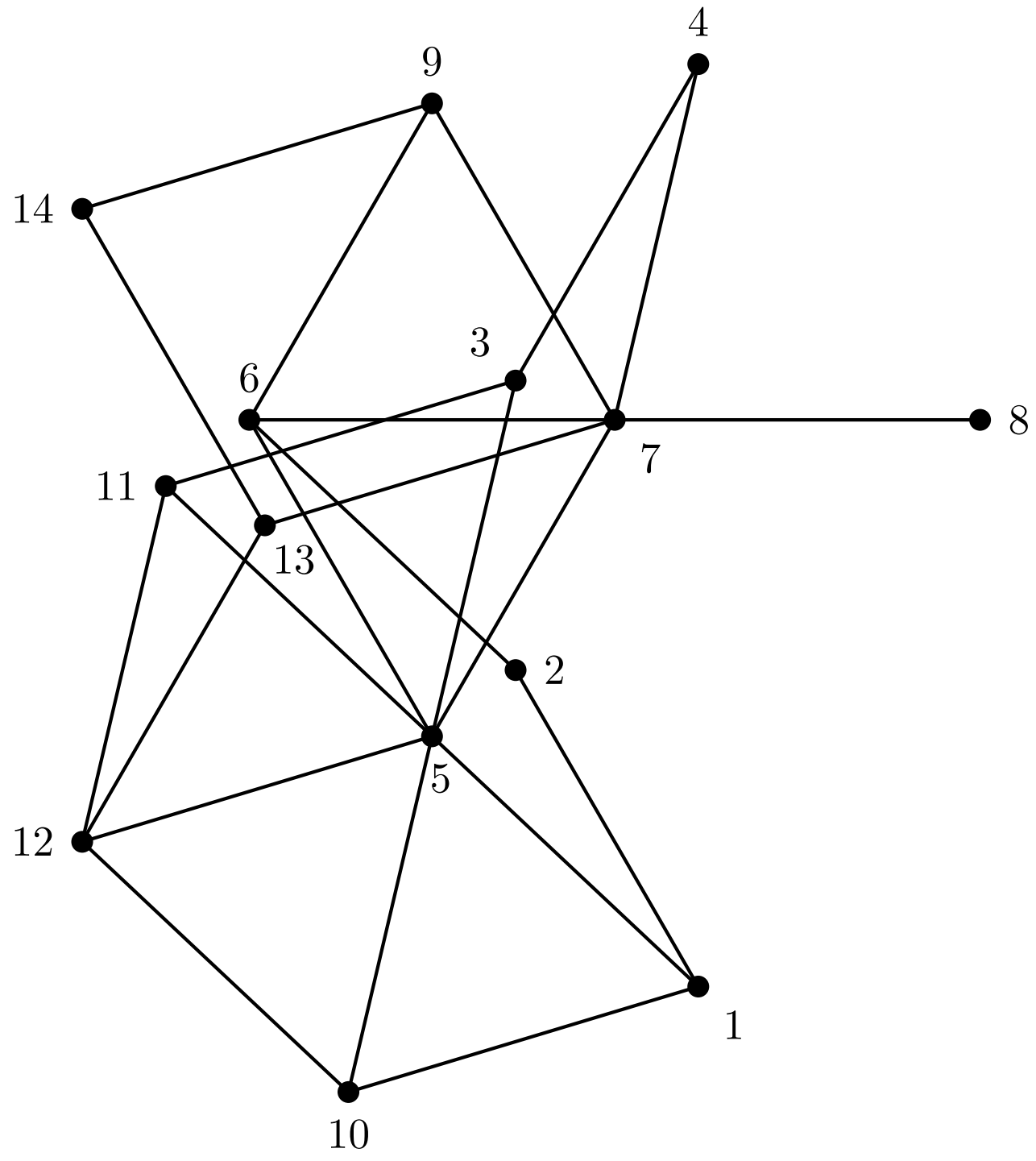
\begin{figure}
  \centering
  \begin{tikzpicture}[
    x=3cm,
    y=3cm,
    vertex/.style={circle, fill=black, inner sep=1.8pt}
  ]
    % Moser-lattice coordinates: (a,b,c,d) represents a+b*omega+c*eta+d*omega*eta,
    % where omega=(1+i*sqrt(3))/2 and eta=(5+i*sqrt(11))/6.
    \foreach \idx/\aa/\bb/\cc/\dd in {
      1/3/3/1/0,
      2/2/4/1/0,
      3/0/3/4/0,
      4/0/4/4/0,
      5/2/2/2/1,
      6/1/3/2/1,
      7/2/3/2/1,
      8/3/3/2/1,
      9/1/4/2/1,
      10/4/1/0/2,
      11/1/1/3/2,
      12/3/0/1/3,
      13/3/1/1/3,
      14/2/2/1/3%
    } {
      \coordinate (v\idx) at
        ({\aa + \bb/2 + 5*\cc/6 + (5-sqrt(33))*\dd/12},
         {sqrt(3)*\bb/2 + sqrt(11)*\cc/6 + (sqrt(11)+5*sqrt(3))*\dd/12});
    }

    % All pairs of vertices at unit distance.
    \foreach \ii/\jj in {
      1/2,
      1/5,
      1/10,
      2/6,
      3/4,
      3/5,
      3/11,
      4/7,
      5/6,
      5/7,
      5/10,
      5/11,
      5/12,
      6/7,
      6/9,
      7/8,
      7/9,
      7/13,
      9/14,
      10/12,
      11/12,
      12/13,
      13/14%
    } {
      \draw[thick] (v\ii) -- (v\jj);
    }

    \node[vertex, label={[font=\small]below right:$1$}] at (v1) {};
    \node[vertex, label={[font=\small]right:$2$}] at (v2) {};
    \node[vertex, label={[font=\small]above left:$3$}] at (v3) {};
    \node[vertex, label={[font=\small]above:$4$}] at (v4) {};
    \node[vertex, label={[font=\small,xshift=2pt]below:$5$}] at (v5) {};
    \node[vertex, label={[font=\small]above:$6$}] at (v6) {};
    \node[vertex, label={[font=\small]below right:$7$}] at (v7) {};
    \node[vertex, label={[font=\small]right:$8$}] at (v8) {};
    \node[vertex, label={[font=\small]above:$9$}] at (v9) {};
    \node[vertex, label={[font=\small]below:$10$}] at (v10) {};
    \node[vertex, label={[font=\small]left:$11$}] at (v11) {};
    \node[vertex, label={[font=\small]left:$12$}] at (v12) {};
    \node[vertex, label={[font=\small,xshift=-4pt,yshift=1pt]below right:$13$}] at (v13) {};
    \node[vertex, label={[font=\small]left:$14$}] at (v14) {};
  \end{tikzpicture}
\caption{The fourteen-point planar probe $X_{14}$, with edges joining pairs of points at unit distance.}
\label{fig_tgfcn_14}
\end{figure}

For $j=1,\dots,14$, let $z_j=a_j+b_j\omega+c_j\eta+d_j\omega\eta$, using $\omega$ and $\eta$ as above.
Figure~\ref{fig_tgfcn_14} depicts this probe, and Table~\ref{tab:tgfcn_14_coordinates} lists the integer coefficients $a_j,b_j,c_j,d_j$.

\begin{table}
\centering
\small
\begin{tabular}{c|rrrr}
\hline
$j$ & $a_j$ & $b_j$ & $c_j$ & $d_j$ \\
\hline
1 & $3$ & $3$ & $1$ & $0$ \\
2 & $2$ & $4$ & $1$ & $0$ \\
3 & $0$ & $3$ & $4$ & $0$ \\
4 & $0$ & $4$ & $4$ & $0$ \\
5 & $2$ & $2$ & $2$ & $1$ \\
6 & $1$ & $3$ & $2$ & $1$ \\
7 & $2$ & $3$ & $2$ & $1$ \\
8 & $3$ & $3$ & $2$ & $1$ \\
9 & $1$ & $4$ & $2$ & $1$ \\
10 & $4$ & $1$ & $0$ & $2$ \\
11 & $1$ & $1$ & $3$ & $2$ \\
12 & $3$ & $0$ & $1$ & $3$ \\
13 & $3$ & $1$ & $1$ & $3$ \\
14 & $2$ & $2$ & $1$ & $3$ \\
\end{tabular}
\vspace{0.3cm}
\caption{Integer coefficients specifying the points of $X_{14}$.
The indices match the labels in Figure~\ref{fig_tgfcn_14}.}
\label{tab:tgfcn_14_coordinates}
\end{table}

\begin{proof}[Proof of Theorem~\ref{thm_plane}]
Every member of $\Pcal_\Delta$ has diameter $1$.
Theorem~\ref{thm:euclidean-probe-bound} then gives $m_\Delta(\R^2)\leq\Gamma_\Delta(X_{14})$.
The linear program gives $\Gamma_\Delta(X_{14})=9/16$.
\end{proof}

\section{Open Problems}
%TODO ezt majd valaki olvassa át hogy tetszik-e. Köszi, Ákos
Several natural questions remain open concerning both extremal densities and the effectiveness of the GFCN method itself.

The most immediate question is whether the density $1/2$ construction described in Section~\ref{sec:plane} is optimal. For unit equilateral triangle avoidance, our upper bound of $\frac{9}{16}$ is still above the conjectured optimal value of $\frac{1}{2}$. Curiously, we were unable to obtain a better upper bound even under the stronger assumption that both a measurable set $A$ and its complement $\R^2\setminus A$ avoid unit equilateral triangles.

The corresponding problem on the sphere also remains open. Williams and Constantin~\cite{WilliamsConstantin2025} conjectured that their scalloped-belt construction is optimal, which would give
\[
m_\perp(\Sphere^2)=1-\frac1{\sqrt2}+\frac{\sqrt2\log 2}{\pi}\approx 0.6049.
\]
We establish the upper bound $7/11\approx 0.6364$. Closing this gap is a natural direction for further work.

A more general question concerns the effectiveness of the GFCN method itself. For a family $\Pcal$ of forbidden configurations in a space $S$, Theorems~\ref{thm:compact-probe-bound} and~\ref{thm:euclidean-probe-bound} establish, under their respective assumptions, that
\[
m_{\Pcal}(S)\le
\inf_{\substack{X\subset S\text{ finite, nonempty}}}
\Gamma_{\Pcal}(X).
\]
For which spaces and forbidden families does equality hold? Equivalently, can the optimal independence density always be approximated arbitrarily closely by bounds obtained from finite probes? For unit-distance avoidance in the plane, the results of~\cite{LiSu2026} and~\cite{Ducz2026} imply equality, since the LP bound for every finite probe is at most its independence ratio. The question remains open in most other cases.

\paragraph{Use of generative AI tools.}

Generative AI tools, including GPT-6 and Claude Fable 5.2, were used
to assist with writing the computer code, drafting the
supplementary documentation, and improving the exposition
and organization of the manuscript. All AI-assisted material was reviewed by the authors, who take full responsibility for the mathematical results, computational verification, and final content of the paper.

\section{Acknowledgements}

D.~C., Á.~D., D.~V., and P.~Zs. were supported by the Hungarian Ministry of Innovation and Technology NRDI Office within the framework of the Artificial Intelligence National Laboratory (RRF-2.3.1-21-2022-00004).
Á.~D. and D.~V. were supported by the grant NKFIH-153165.
D.~C. and P.~M. were supported by the grant NKFIH-154121.
We thank András Szabó for his help in preparing the supplementary material, and Gergely Ambrus, Máté Matolcsi, and Kenneth Moore for helpful discussions and valuable insights.

\appendix

\section{Computational search for the probes}\label{app:search}
% TODO: Valakinek ellenőriznie kell ezt a fejezetet, aki érti!
% Ákos: nekem jónak tűnik, amennyiben a projektivitással minden oké, erre volt egy korábbi TODO

The probes $X_{14}$ and $X_{15}$ were obtained by computer-assisted searches over finite geometric configurations, following the general approach of~\cite{AmbrusCsiszarikMatolcsiVargaZsamboki2024,MatolcsiRuzsaVargaZsamboki2025,CzifraEtAl2026}.
The search software generates candidate configurations by geometric extensions and ranks them using numerical solutions of their GFCN linear programs.
Both beam search and Monte Carlo tree search are implemented for exploring the resulting space of configurations.

For the planar problem, points are represented by integer quadruples $(a,b,c,d)$ encoding $a+b\omega+c\eta+d\omega\eta$, where $\omega$ and $\eta$ are defined in Section~\ref{sec:plane}.
Candidate extensions include adjoining points at unit distance from existing points, completing equilateral triangles, and completing parallelograms.
Configurations are canonicalized under translations and symmetries of the Moser lattice to avoid duplicate candidates.
The fourteen integer quadruples for the resulting probe $X_{14}$ are listed in Table~\ref{tab:tgfcn_14_coordinates}.

For the spherical problem, the search works with projective directions, identifying a unit vector $u$ with $-u$.
The extension operation takes two nonparallel representatives $u,v$ and adjoins the projective direction represented by $u\times v$, which is orthogonal to both.
The same cross-product extension rule is used in the spherical search described in~\cite{CzifraEtAl2026}.
The exact coordinates of $X_{15}$ are listed in Table~\ref{tab:vertices-projective}.

For each candidate $X$, its graph has edges between pairs at unit distance in the planar problem or between orthogonal pairs in the spherical problem.
The triangle-free vertex subsets are precisely the admissible observation patterns, and the linear program also imposes the congruence constraints of Sections~\ref{sec:compact} and~\ref{sec:euclidean}.
The search evaluates the LP numerically, using a smaller value of $\Gamma_\Delta(X)$ or $\Gamma_\perp(X)$ as a better score.
For the final probes, $X_{14}$ has $23$ unit-distance pairs and $6$ unit equilateral triangles, while $X_{15}$ has $25$ orthogonal pairs and $6$ orthonormal triples.
These counts can be checked directly from the exact coordinates, whereas the LP congruence constraints also use geometric information not contained in the graphs.
The exact LP values $\Gamma_\Delta(X_{14})=9/16$ and $\Gamma_\perp(X_{15})=7/11$ are verified by rational primal and dual certificates independently of the numerical search.

\end{document}